\documentclass[12pt]{article}\usepackage{amsfonts,amssymb}

\usepackage{amsmath,bbm}
\usepackage{theorem}          
\usepackage{latexsym} 
\usepackage{float}                   
\usepackage{times}
\usepackage{cite}

\newtheorem{tht}{Theorem}[section]

\newtheorem{thl}[tht]{Lemma}

\newcommand{\hs}{\hspace{-1pt}}
\newcommand{\anf}{\raisebox{0.2ex}{\scriptsize$\triangleleft$}}
\newcommand{\ang}{\raisebox{0.2ex}{\scriptsize$\triangleright$}}

\newcommand{\mn}{\medskip}    

\newcommand{\sn}{\smallskip\noindent}

\newcommand{\rti}{\,{\scriptstyle\rtimes}\,} %

\newcommand{\im}{\mathrm{i}}
\newcommand{\tr}{\mathrm{Tr}}
\newcommand{\dd}{\mathrm{d}}  
\newcommand{\su}{\mathrm{su}}
\newcommand{\OSU}{\mathrm{SU}}
\newcommand{\cD}{{\mathcal{D}}}

\newcommand{\Hh}{{\mathcal{H}}} 
\newcommand{\cO}{{\mathcal{O}}}

\newcommand{\cU}{{\mathcal{U}}}   

\newcommand{\cI}{\mathcal{I}}
\newcommand{\cX}{{\mathcal{X}}}
\newcommand{\cA}{{\mathcal{A}}}

\newcommand{\cC}{{\mathcal{C}}}
\newcommand{\cF}{{\mathcal{F}}}

\newcommand{\cL}{\mathcal{L}}

\newcommand{\dZ}{\mathbb{Z}}

\newcommand{\dC}{\mathbb{C}}
\newcommand{\dN}{\mathbb{N}}

\newcommand{\gk}{\mathfrak{k}}

\newcommand{\gj}{{\mathfrak{j}}}
\newcommand{\gl}{\mathfrak{l}}

\newcommand{\inv}{\mathrm{inv}}
\newcommand{\Lin}{{\mathrm{Lin}}}
\newcommand{\id}{{\mathrm{id}}}

\begin{document}

\date{\small{Fakult\"at f\"ur Mathematik und 
Informatik\\ Universit\"at Leipzig, 
Augustusplatz 10, 04109 Leipzig, Germany\\ 
E-mail: schmuedg@mathematik.uni-leipzig.de / 
wagner@mathematik.uni-leipzig.de}
}

\title{Examples of twisted cyclic cocycles from covariant differential calculi}

\author{Konrad Schm\"udgen and Elmar Wagner}

\maketitle

\renewcommand{\theenumi}{\roman{enumi}}
\begin{abstract}\noindent
For two covariant differential $\ast$-calculi, the twisted cyclic 
cocycle associated with the volume form is represented in terms of 
commutators $[\cF,\rho(x)]$ for some self-adjoint operator $\cF$ 
and some $\ast$-re\-pre\-sen\-tation $\rho$ of the underlying 
$\ast$-al\-ge\-bra.
\end{abstract}

\setcounter{section}{-1}
\section{Introduction}
Twisted cyclic cocycles arise under certain assumptions from covariant 
differential calculi on quantum groups or quantum spaces \cite{KMT}. 
Let $\Gamma^\wedge = \oplus_k\Gamma^{\wedge k}$ be a differential calculus 
on an algebra $\cX$. We assume the existence of a volume 
form $\omega \in\Gamma^{\wedge n}$ (that is, for each $n$-form 
$\eta\in\Gamma^{\wedge n}$, there is a unique element $\pi(\eta)\in\cX$ 
such that $\eta = \pi(\eta)\omega$),
of an algebra automorphism $\sigma_1$ of $\cX$ for which    
$\omega x= \sigma_1(x)\omega$,  
 and  of a twisted trace $h$ on $\cX$ 
(that is, there is an algebra automorphism $\sigma_2$ of $\cX$ such 
that $h(xy) = h(\sigma_2(y)x)$ for all $x, y\in\cX)$. The twisted cyclic 
cocycle $\tau_{\omega , h}$ on $\cX$ associated with $\omega $ and $h$ 
is defined by
$$
\tau_{\omega , h} (x_0,  x_1,\dots, x_n) 
= h(\pi(x_0 \dd x_1\wedge\dots\wedge \dd x_n)),\quad  x_0, x_1,\dots,  
x_n\in \cX.
$$
Precise definitions and assumptions are given below.

The purpose of this paper is to investigate two examples of such 
twisted cyclic cocycles $\tau_{\omega ,  h}$. The first example, discussed 
in Section 2, is the $\ast$-al\-ge\-bra $\cX_0$ 
with generators $z$, $z^\ast$ and defining relation 
$z^\ast z\hs -\hs q^2 zz^\ast \hs=\hs (1\hs-\hs q^2)\alpha$, 
where $0\hs<\hs q\hs <\hs 1$ and 
$\alpha = 0,  1,  -1$. We study the twisted cyclic cocycle associated 
with the dis\-tin\-guished differential calculus on this $\ast$-al\-ge\-bra. 
The second example, treated in Section 3, concerns left-covariant 
differential calculi on Hopf algebras. First we develop the general 
framework of representing differential forms in terms of commutators 
and then we carry out the details for Woronowicz' 3D-calculus on 
the quantum group $\OSU_q(2)$. Our main technical tool is to construct 
an appropriate commutator representation $\dd x\cong \im [\cF,  \rho(x)]$ of 
the corresponding first order $\ast$-cal\-cu\-lus, where $\rho$ is 
a $\ast$-re\-pre\-sen\-tation of $\cX$ and $\cF$ 
is a self-adjoint operator 
on a Hilbert space. In both examples, the cocycle $\tau_{\omega ,h}$ is then 
described in the form
\begin{equation}\label{tracecom}
\tau_{\omega ,h}(x_0,  x_1,\dots,  x_n) 
= \tr\,  A \gamma_q \rho (x_0) [\cF,  \rho (x_1)]\cdots  [\cF,  \rho (x_n)],  
\end{equation}
where $A$ is a certain density operator and $\gamma_q$ 
is a grading operator (see Theorems 2.1 and 3.3 below for details).
%
\section{Twisted cyclic cohomology}\label{S1}
%
Suppose that $\cX$ is a complex algebra and $\sigma$ is 
an algebra automorphism of $\cX$.
Let $\varphi$ be an $(n+1)$-linear form on $\cX$. 
The $\sigma$-twisted coboundary operator $b_\sigma$ 
and the $\sigma$-twisted cyclicity operator $\lambda_\sigma$ 
on $\cX$ are defined by
\begin{align*}
(b_\sigma\varphi)(x_0,\dots,  x_n) 
&= \sum^{n-1}_{j=0} (-1)^j \varphi (x_0,\dots,  x_j x_{j+1},\dots,  x_n)\\
&\qquad\qquad\qquad  
+(-1)^n \varphi (\sigma(x_n)x_0,  x_1,\dots, x_{n-1}),\\[6pt]
(\lambda_\sigma \varphi)(x_0,\dots,  x_n) 
&= (-1)^{n} \varphi (\sigma(x_n),  x_0,\dots,  x_{n-1}),
\end{align*}
where $x_0,\dots,  x_n\in\cX$. An $(n+1)$-form $\varphi$ 
is called a {\it $\sigma$-twisted cyclic $n$-cocycle} 
if $b_\sigma\varphi = 0 $ and $\lambda_\sigma\varphi =\varphi$.

Twisted cyclic $n$-cycles occur in the study 
of differential calculi on algebras (see \cite{KMT} and \cite{KR}). 
Let $\Gamma^\wedge = \oplus^\infty_{k=0} \Gamma^{\wedge k}$ 
be a differential calculus on $\cX$ with differentiation 
$\dd : \Gamma^{\wedge k}\rightarrow\Gamma^{\wedge(k+1)}$. 
We assume that there is an $n$-form $\omega \in\Gamma^{\wedge n}$ 
such that for each $\eta\in\Gamma^{\wedge n}$ there exists 
a unique element $\pi(\eta)\in\cX$ such that $\eta = \pi(\eta)\omega$. 
Moreover, we assume that there is an algebra automorphism $\sigma_1$ 
of $\cX$ such that
\begin{equation}\label{omegaaut}
\omega x = \sigma_1 (x) \omega,\quad x\in \cX.
\end{equation}
Suppose further that $h$ is a linear functional on $\cX$ and $\sigma_2$ 
is an algebra automorphism of $\cX$ satisfying
\begin{equation}\label{haut}
h(x_1x_2) = h(\sigma_2(x_2)x_1),\quad x_1,  x_2\in\cX,
\end{equation}
\begin{equation}\label{hinv}
h(\pi(\dd x_1\wedge\dots\wedge \dd x_n)) = 0,\quad  x_1,\dots,  x_n\in\cX.
\end{equation}
Define $\sigma = \sigma_2\circ\sigma_1$ and 
\begin{equation}\label{tau}
\tau_{\omega,h}(x_0,  x_1,\dots, x_n) 
:= h(\pi(x_0 \dd x_1\wedge\dots\wedge x_n)),\quad x_0,\dots,  x_n\in\cX.
\end{equation}
Then, using the Leibniz rule and assumption (\ref{omegaaut}) 
and (\ref{haut}), it is not difficult to check that 
$b_\sigma\tau_{\omega,h} = 0$. Using in addition (\ref{hinv}), 
we obtain $\lambda_\sigma\tau_{\omega,h} =\tau_{\omega,h}$. 
(Note that (\ref{hinv}) 
is crucial for the proof of the relation 
$\lambda_\sigma\tau_{\omega,h} =\tau_{\omega,h}$.) 
Thus, $\tau_{\omega,h}$ is a $\sigma$-twisted cyclic $n$-cycle on $\cX$. 
We call $\tau_{\omega,h}$ the twisted cyclic $n$-cycle associated 
with the $n$-form $\omega$ of the differential calculus $\Gamma^\wedge$ 
and the functional $h$.

Let us briefly discuss the preceding assumptions on the calculus. 
Assumptions (\ref{haut}) and (\ref{hinv}) are fulfilled in most 
interesting cases, for instance, if $\Gamma^\wedge$ is a left or right 
covariant differential calculus on a compact quantum group 
Hopf algebra $\cX$ and $h$ is the Haar state of $\cX$. 
The existence of a form $\omega\in\Gamma^{\wedge n}$ as above is 
satisfied for many, but not all, covariant differential calculi 
on Hopf algebras. It is satisfied for the standard bicovariant 
differential calculi on the quantum groups $\mathrm{GL}_q(k)$ and 
$\mathrm{SL}_q(k)$ 
(with $n = k^2$ and $\sigma_1 = \id$, see \cite{S3}). 
The paper \cite{H} contains a list of left  covariant differential calculi 
on $\cO(\mathrm{SL}_q(2))$ satisfying the assumptions with $n=3$. 
%
%
\section{Quantum disc, complex quantum plane, and ster\-eo\-gra\-phic 
projection of the Podles' sphere} 
                                                           \label{S2}
In this section, let $\cX_0$ be the unital $\ast$-al\-ge\-bra 
with generators $z$, $z^\ast$ and defining relation 
\begin{equation}\label{zrel}
z^\ast z - q^2 zz^\ast = \alpha(1-q^2),
\end{equation}
where $q$ and $\alpha$ are fixed real numbers such that $0<q<1$ 
and $\alpha\in\{0,1,-1\}$. There is a distinguished first order 
differential $\ast$-calculus $(\Gamma_0,\dd)$ 
on $\cX_0$ with $\cX_0$-bi\-mo\-dule 
structure given by 
\begin{equation}\label{dzrel}
\dd z\, z = q^2 z \dd z,\ \ \dd z\, z^\ast = q^{-2} z^\ast \dd z, \ \ 
\dd z^\ast \, z = q^2 z \dd z^\ast, \ \ 
\dd z^\ast \, z^\ast = q^{-2} z^\ast \dd z^\ast.
\end{equation}
These simple relations have been found in \cite{S2}.

All three cases $\alpha = 0,  1,  -1$ give interesting quantum spaces. 
In the case $\alpha = 0$, we get the quantum complex plane 
on which the Hopf $\ast$-al\-ge\-bra $\cU_q(\mathrm{e}_2)$ acts. 
For $\alpha = 1$, it gives the quantum disc 
algebra \cite{KL,SSV} which  
is a $\cU_q(\su_{1,1})$-mo\-dule $\ast$-al\-ge\-bra. 
For $\alpha = -1$, the $\ast$-al\-ge\-bra $\cX_0$ is a left 
$\cU_q(\su_2)$-mo\-dule $\ast$-al\-ge\-bra which can be interpreted 
as the stereographic projection of the standard Podles' quantum sphere. 
To be more precise, the left action of  $\cU_q(\su_2)$ 
on the Hopf $\ast$-al\-ge\-bra $\cO(\OSU_q(2))$ of the quantum $\OSU(2)$ group 
extends to an action on the Ore extension $\hat{\cO}(\OSU_q(2))$ 
with respect to the Ore set $\{b^r c^s\,;\,r,  s\in\dN_0\}$. 
Then the left $\cU_q(\su_2)$-mo\-dule $\ast$-al\-ge\-bra $\cX_0$ 
can be identified with the $\ast$-sub\-al\-ge\-bra of $\hat{\cO}(\OSU_q(2))$ 
generated by the element $z: = ac^{-1}$. In all three cases, the 
differential calculus $(\Gamma_0,\dd)$ on $\cX_0$ is covariant 
with respect to the action of the corresponding Hopf $\ast$-al\-ge\-bras.

Next we pass to an appropriate Hilbert space representation 
of the $\ast$-al\-ge\-bra $\cX_0$. It acts on an orthonormal basis 
$\{e_n\,;\,n\in\dZ\}$ of $\Hh = l^2(\dZ)$ 
and  $\{e_n\,;\,n\in\dN_0\}$ of $\Hh = l^2(\dN_0)$ for 
$\alpha = 0$ and $\alpha = \pm 1$,
respectively, by the following formulas:
\begin{align*}
\alpha &= 0 : & ze_n &=  q^{2(n+1)} e_{n+1},& z^\ast e_n &= q^{2n} e_{n-1},\\
\alpha &= 1 : & ze_n &=  (1-q^{2(n+1)})^{1/2} e_{n+1}, &
z^\ast e_n &= (1-q^{2n})^{1/2} e_{n-1},\\
\alpha &= -1 : &  ze_n &=  (q^{-2n}-1)^{1/2} e_{n-1}, & 
z^\ast e_n &= (q^{-2(n+1)}-1)^{1/2} e_{n+1}.
\end{align*}
Let $\cX$ be the $\ast$-al\-ge\-bra $\cL^+(\cD)$ of all linear operators $T$ 
on $\cD = \Lin\{e_n\,;\,n\in\dZ\}$ resp.\ $\cD = \Lin\{e_n\,;\,n\in\dN_0\}$
such that $T\cD\subseteq\cD$ and 
$T^\ast\cD\subseteq\cD$. We identify the $\ast$-al\-ge\-bra $\cX_0$ 
with the corresponding $\ast$-sub\-al\-ge\-bra of $\cX$.

Set $\beta = 1$ for $\alpha = 1$ and $\beta = -1$ for $\alpha = 0,  -1$. 
Put $y:=\beta(\alpha-zz^\ast)$. 
Note that $y^{-1}\in \cX$. 
Let $\cX_c$ be the set of all elements $x\in\cX$ such that 
$y^kxy^l$ is bounded for all $k,l\in\dZ$. 
Then $\cX_c$ is a $\ast$-sub\-al\-ge\-bra of $\cX$ and, moreover, 
a $\cX_0$-bimodule. Notice that the closures of all operators $y^kxy^l$, 
$x\in\cX_c$, are of trace class. 
We define an algebra automorphism 
$\sigma$ of $\cX$ by $\sigma(x) = yxy^{-1}$, $x\in\cX$. Since $y\ge 0$,
\begin{equation}\label{hdef}
h(x) := \tr\, y^{-1}x,\quad x\in\cX_c,
\end{equation}
is a positive linear functional on $\cX_c$ which obviously satisfies 
$$
h(x_1x_2) = h(\sigma(x_2)x_1),\quad x_1, x_2\in\cX_c.
$$
It can be shown that the action of the Hopf $\ast$-al\-ge\-bra 
$\cU_q(\mathrm{e}_2)$, 
$\cU_q(\su_{1,1})$ resp.\ $\cU_q(\su_2)$ on $\cX_0$ extends to $\cX$ such 
that $\cX$ and $\cX_c$ are left module $\ast$-al\-ge\-bras and $h$ is 
an invariant functional on $\cX_c$ (cf.\ \cite{KW}). 
We do not carry out the details, 
because we will not need this fact in what follows.

Next we extend the first order differential calculus $(\Gamma_0,\dd)$ 
on $\cX_0$ to the larger $\ast$-al\-ge\-bra $\cX$ by using a commutator 
representation (see \cite{S1}). Define an operator $\cF$ 
on $\Hh\oplus \Hh$ and a $\ast$-homomorphism $\rho$ of $\cX$ 
into the $\ast$-al\-ge\-bra $\cL^+(\cD\oplus\cD)$ by
$$
\cF = (1-q^2)^{-1}\left( \begin{matrix} 0 &z\\ 
z^\ast &0 \end{matrix} \right), \quad \rho(x) = \left( \begin{matrix} x &0\\ 
0 &x \end{matrix} \right), \quad x\in\cX.
$$
Then there is an injective linear map 
$$
\Gamma_0\ni x_1\dd x_2\rightarrow 
\rho(x_1)[\im \cF, \rho(x_2)] \in \cL^+(\cD\oplus\cD),\quad  x_1, x_2\in\cX_0.
$$ 
For simplicity, 
we identify $x_1 \dd x_2$ with $\rho(x_1)[\im \cF, \rho(x_2)]$, 
so the differentiation of the calculus $(\Gamma_0,\dd)$ on $\cX_0$ 
is given by the commutator with the operator $\im \cF$. 
Defining $x_1 \dd x_2 := \rho(x_1) [\im \cF, \rho(x_2)]$ for $x_1, x_2\in\cX$, 
the calculus $(\Gamma_0,\dd)$ is extended to a first order 
$\ast$-calculus $(\Gamma,\dd)$ on the $\ast$-al\-ge\-bra $\cX$.
Since $[z^\ast, z] = (1-q^2)(\alpha-zz^\ast)=(1-q^2)\beta y$ 
by (\ref{zrel}), we have 
$$
\dd z  = \left( \begin{matrix} 0&0\\ \im \beta y &0 \end{matrix} \right),
\qquad \dd z^\ast 
= \left( \begin{matrix} 0 &-\im \beta y\\ 0 &0 \end{matrix} 
\right),
$$
so that
\begin{equation}\label{dzx}
\dd z\, x = \sigma(x) \dd z,\quad 
\dd z^\ast\, x = \sigma(x) \dd z^\ast, \quad x\in\cX.
\end{equation}
Since $y$ is invertible, 
$\{\dd z, \dd z^\ast\}$ is a left module basis of $\Gamma$.
Thus, for any $x\in\cX$, 
there are uniquely determined elements $\partial_z(x)$, 
$\partial_{z^\ast}(x)\in\cX$ such that 
$$
\dd x = \partial_z(x) \dd z +\partial_{z^\ast} (x) \dd z^\ast.
$$
Comparing the latter and the expressions for $\dd z$ and $\dd z^\ast$, we get 
\begin{equation}\label{deltarel}
\partial_z(x) = (1-q^2)^{-1}\beta [z^\ast, x] y^{-1},\quad 
\partial_{z^\ast} (x) = -(1-q^2)^{-1} \beta [z, x]y^{-1}.
\end{equation}
Let us describe the higher order differential calculus
$\Gamma^\wedge = \oplus_k \Gamma^{\wedge k}$ on $\cX$ associated 
with the first order calculus $(\Gamma,\dd)$. The non-zero 2-form 
$\omega := q^{-2} y^{-2} \dd z^\ast\wedge \dd z$ is a left module basis of 
$\Gamma^{\wedge 2}$ and we have $\Gamma^{\wedge k}  = \{0\}$ if $k\ge 3$. 
(Incidentally, $\omega$ is invariant with respect to the action 
of Hopf $\ast$-al\-ge\-bras mentioned above.) From (\ref{dzx}) and 
the definition of $\sigma$,  
we get $\omega x = x\omega$ for all $x\in\cX$, that is, 
condition (\ref{omegaaut}) is satisfied with $\sigma_1$ 
being the identity. Moreover, since 
$\dd z\wedge \dd z =\dd z^\ast\wedge \dd z^\ast = 0$ 
and $\dd z^\ast \wedge \dd z = -q^2 \dd z \wedge \dd z^\ast$ 
in $\Gamma^{\wedge 2}$ 
by  (\ref{dzrel}), it follows that
$$
x_0 \dd x_1 \wedge \dd x_2 
=  x_0(q^2 \partial _{z^\ast} (x_1) \sigma (\partial_z (x_2)) 
- \partial_z (x_1) \sigma (\partial_{z^\ast} (x_2))) y^{2} \omega
$$
for $x_0$, $x_1$, $x_2\in\cX$. Therefore, using (\ref{tau}), (\ref{hdef}) 
and (\ref{deltarel}), we compute
\begin{align}\label{tautrace}
\tau_{\omega,h} (x_0,  x_1, x_2)
&= h(x_0 (q^2 \partial_{z^\ast} (x_1)\sigma (\partial_z (x_2)) 
- \partial_z(x_1)\sigma (\partial_{z^\ast} (x_2)))y^2) \nonumber \\
&= (1-q^2)^{-2} \tr\,  y^{-1} x_0(-q^2 [z, x_1][z^\ast, x_2] 
+ [z^\ast, x_1][z, x_2])\qquad \nonumber \\
&= \tr\,  \gamma_q \rho (y^{-1}) \rho (x_0) [\cF, \rho (x_1)][\cF, \rho(x_2)]
\end{align}
for $x_0$, $x_1$, $x_2\in\cX_c$, where $\gamma_q$ denotes  
the ``grading operator" 
$$
\gamma_q = \left( \begin{matrix} q^2 &0\\ 0 &-1\end{matrix} \right).
$$
Finally, we have to show that condition (\ref{hinv}) is satisfied. 
This could be done by using the invariance of $\omega$ under the Hopf 
algebra action, but we prefer to check it directly. 
Using the commutation relations $y^{-1} z = q^{-2} zy^{-1}$, 
$y^{-1} z^\ast = q^2 z^\ast y^{-1}$ and (\ref{zrel}), we compute
\begin{align}\label{verhinv}
&\tr\,  y^{-1} ([z^\ast, x_1][z, x_2]-q^2[z, x_1][z^\ast, x_2])
\nonumber\\
&\quad =\tr\,  y^{-1} (-x_1 z^\ast z x_2+q^2x_1 zz^\ast x_2)-
\tr\,  y^{-1} (z^\ast x_1 x_2 z - q^2 zx_1x_2z^\ast)\nonumber\\
&\quad\qquad +\tr\,  y^{-1} (z^\ast x_1zx_2 - q^2 x_1 zx_2z^\ast) 
+ \tr\,  y^{-1} (x_1z^\ast x_2 z - q^2zx_1z^\ast x_2)\nonumber\\
&\quad =\tr\,  y^{-1} (-x_1 x_2 \alpha (1-q^2)) 
- \tr\,  (q^{2} z^\ast y^{-1} x_1 x_2 z - z y^{-1} x_1 x_2 z^\ast)\nonumber\\
&\quad\qquad 
+\tr\, (q^2 z^\ast y^{-1} x_1 z x_2 - q^2 y^{-1}x_1 z x_2 z^\ast) 
+\tr\, ( y^{-1} x_1 z^\ast x_2 z -zy^{-1}x_1 z^\ast x_2) \nonumber\\
&\quad = 0
\end{align}
for $x_1,x_2\in\cX_c$. 
The last equality follows from the trace property since 
the closures of the operators  $y^{-1}x_1$, $zy^{-1}x_1$, $z^\ast y^{-1}x_1$
are of trace class and $x_2$,  $x_2z$, $zx_2$, $x_2 z^\ast$, 
$z^\ast x_2$, $z^\ast x_2 z$, $z x_2 z^\ast$ are bounded by the definition 
of $\cX_c$. As
$$ h(\pi(\dd x_1 \dd x_2)) = \tau_{\omega,h}(1, x_1, x_2) 
=  \tr\,  \gamma_q\rho(y^{-1}) [\cF, \rho(x_1)] [\cF, \rho(x_2)] = 0
$$
by (\ref{tautrace}) and (\ref{verhinv}), condition (\ref{hinv}) holds. 
Therefore, $\tau_{\omega,h}$ is a $\sigma$-twisted cyclic 2-cocycle.

We summarize the preceding in the following theorem.
\begin{tht} \label{T2/2}
In the above notation, $\tau_{\omega,h}$ is a $\sigma$-twisted 
cyclic 2-cycle on the $\ast$-al\-ge\-bra $\cX_c$ and we have
$$
\tau_{\omega,h} (x_0, x_1, x_2) 
= \tr\, \gamma_q\rho(y^{-1}) \rho(x_0) [\cF, \rho (x_1)][\cF, \rho(x_2)],
\quad x_0, x_1, x_2\in\cX_c.
$$
\end{tht}
%
\section{The 3D-Calculus on the Quantum Group ${\bf \OSU_q(2)}$} \label{S3}
%
In this section, we use some facts on covariant differential calculi 
on Hopf algebras which can be found e.g.\  in \cite[Chapter 14]{KS}. 
Suppose that $(\Gamma, \dd)$ is a finite dimensional 
left-covariant first order differential calculus 
on a Hopf algebra $\cA$ and $\Gamma^\wedge  = \oplus_j \Gamma^{\wedge j}$ 
is a differential calculus on $\cA$  such that $\Gamma^{\wedge 1} = \Gamma$. 
As in Section 1, we assume that there are a natural number $n$ 
and a $n$-form $\omega\in\Gamma^{\wedge n}$ such that 
for each $\xi\in\Gamma^{\wedge n}$ there exists a unique element 
$\pi(\xi)\in\cA$ such that $\xi = \pi(\xi)\omega$. 
Our first aim is to express $\pi(x_0 \dd x_1\wedge\cdots \wedge \dd x_n)$ 
in terms of commutators in the cross product algebra $\cA\rti \cA^0$.

Let $\{\omega_1,\dots, \omega_m\}$ be a basis 
of the vector space $\Gamma_{\inv}$ of left-in\-var\-iant elements 
of $\Gamma$. Then there exist functionals $X_k$, $f^k_l$, 
$k,l = 1,\dots, m$, of the Hopf dual $\cA^\circ$ of $\cA$ such that
\begin{align}\label{dxfcom}
\dd x &= \sum_k(X_k\ang x) \omega_k,\quad \omega_k x 
= \sum_j (f^k_j\ang x) \omega_j\\
\label{xfco}
\Delta(X_k) &= \varepsilon \otimes X_k+\sum_jX_j\otimes f^j_k,\quad 
\Delta(f^k_l) = \sum_j f^k_j \otimes f^j_l, 
\end{align}
where $f\ang x := \langle f, x_{(2)}\rangle x_{(1)}$ is 
the left action of $f\in\cA^\circ$ on $x\in\cA$. Recall that 
the cross product algebra $\cA\rti\cA^\circ$ is the algebra 
generated by the two subalgebras $\cA$ and $\cA^\circ$ 
with respect to the cross relation
\begin{equation}\label{cross}
fx = (f_{(1)}\ang x) f_{(2)} 
\equiv \langle f_{(1)},x_{(2)}\rangle x_{(1)} f_{(2)}, \quad 
x\in\cA,\ \, f\in\cA^\circ.
\end{equation}

We shall use the notation $(\gk, l) := (k_1,\dots, k_r, l)$, 
$\omega_\gk := \omega_{k_1} \omega_{k_2}\cdots \omega_{k_r}$, 
$f^\gk_\gj := f^{k_1}_{j_1}\cdots f^{k_r}_{j_r}$ for multi-indices 
$\gk = (k_1,\dots, k_r)$ and $\gj = (j_1\dots, j_r)$, 
where $k_i,j_i, l\in\{1,\dots, m\}$.
Let $x_0, x_1,\dots, x_n\in\cA$. From (\ref{dxfcom}), we deduce 
\begin{align}                                           \label{xdrel}
&x_0 \dd x_1\wedge\cdots \wedge \dd x_n\\ \nonumber
&\quad
=\sum_{\gk_i,l_i} x_0(X_{l_1} \ang x_1)(f_{\gk_1}^{l_1} X_{l_2} \ang x_2)
(f^{(\gk_1,l_2)}_{\gk_2} X_{l_3} \ang x_3)
\cdots (f^{(\gk_{n-2},l_{n-1})}_{\gk_{n-1}} X_{l_n} \ang x_n) 
\omega_{(\gk_{n-1},l_n)},
\end{align}
where the summation is over all multi-indices 
$\gk_i = (k_{i1},\dots, k_{ii})$ and numbers $l_i\in \{1,\dots, m\}$.
Furthermore, Equations (\ref{xfco}) and (\ref{cross}) give 
\begin{equation}\label{xfcomma}
[X_k, x] = \sum_j (X_j\ang x) f^j_k,\quad f^k_j x  
= \sum_l(f^k_l \ang x)f^l_j,\quad x\in\cA. 
\end{equation}
Using these relations, a straightforward induction argument shows that 
\begin{align*}
&[X_{j_1}, x_1] [X_{j_2}, x_2]\cdots [X_{j_n}, x_n]\\
&\qquad =\sum_{\gk_i,l_i} (X_{l_1} \ang x_1) 
(f^{l_1}_{\gk_1} X_{l_2} \ang x_2) 
(f^{(\gk_1, l_2)}_{\gk_2} X_{l_3} \ang x_3)
\cdots (f^{(\gk_{n-2},l_{n-1})}_{\gk_{n-1}} X_{l_n} \ang x_n) 
f^{(\gk_{n-1},l_n)}_{(j_1,\dots,j_n)}
\end{align*}
for $j_1,\dots, j_n\in\{1,\dots, m\}$. Since 
$\sum_{\gj}f^\gk_\gj S(f^\gj_\gl)=\varepsilon(f^\gk_\gl)=\delta^\gk_\gl$, 
we get 
\begin{align}             \label{comrel}
&\sum_{\gj=(j_1,\ldots,j_n)} 
[X_{j_1}, x_1]\cdots [X_{j_n}, x_n] S(f^\gj_{(\gk,l)})\\ \nonumber 
&\quad\quad=\sum_{\gk_i,l_i} (X_{l_1}\ang x_1) 
(f^{l_1}_{\gk_1} X_{l_2}\ang x_2)
\cdots (f^{(\gk_{n-3},l_{n-2})}_{\gk_{n-2}}                 
X_{l_n-1} \ang x_{n-1}) (f^{(\gk_{n-2},l_{n-1})}_{\gk} X_l \ang x_n).
\end{align}
Combining (\ref{xdrel}) and (\ref{comrel}) 
proves the assertion of the following lemma.
\begin{thl} \label{T3/1}
For $x_0, x_1,\dots, x_n\in\cA$, we have 
$$
\pi(x_0 \dd x_1\wedge \cdots  \wedge \dd x_n) 
= \sum_{j_i,k_i} x_0 [X_{j_1}, x_1] \cdots  [X_{j_n}, x_n] 
S(f^{j_1}_{k_1}\cdots f^{j_n}_{k_n}) \pi 
(\omega_{k_1}\cdots \omega_{k_n}).
$$
\end{thl}
As in Section \ref{S2}, we want to express 
the twisted cyclic cocycle $\tau_{\omega,h}$ in terms 
of commutators $[\cF, \rho(x)]$, where $\cF$ is a self-adjoint 
operator and $\rho$ is a representation of $\cA$ on a Hilbert space. 
It seems reasonable  
to look for an operator $\cF$ of the form 
$\cF = \sum_j X_j\otimes \eta_j$, where $\eta_1,\dots, \eta_m$ 
are appropriate matrices. We shall carry out this 
for Woronowicz' 3D-calculus \cite{W} on the quantum group $\OSU_q(2)$.
Let $\cA$ be the Hopf $\ast$-al\-ge\-bra $\cO(\OSU_q(2))$ 
with usual generators $a$, $b$, $c$, $d$ and let $\cU_q(\su_2)$ 
be the Hopf $\ast$-al\-ge\-bra with generators $E$, $F$, $K$, $K^{-1}$, 
relations
\begin{align*} 
&KK^{-1} = K^{-1}K = 1,\quad  KE = qEK, \quad FK = qKF, \\
&EF-FE = (K^2-K^{-2})/(q-q^{-1}),
\end{align*}
involution $E^\ast = F$, $K^\ast = K$, comultiplication 
$$
\Delta(E) = E\otimes K+K^{-1} \otimes E,\ \ \Delta(F) 
= F\otimes K +K^{-1} \otimes F,\ \  \Delta(K) = K\otimes K, 
$$
counit $\varepsilon (E) = \varepsilon (F)= \varepsilon (K-1)=0$ and 
antipode  $S(K)=K^{-1}$, $S(E)=-qE$, $S(F)=-q^{-1}F$. 
There is a non-degenerate dual pairing of these Hopf $\ast$-al\-ge\-bras 
given on ge\-ne\-ra\-tors by 
$\langle K^{\pm 1}, d\rangle = \langle K^{\mp 1}, a\rangle = q^{\pm 1/2}$, 
$\langle E, c\rangle  = \langle F, b\rangle  = 1$ and zero otherwise. 
For the Haar state $h$ on $\cO(\OSU_q(2))$, we have
$$ 
h(x_1x_2) = h(\sigma_2(x_2)x_1)\ \  \mbox{with}\ \ 
\sigma_2(x) = K^{-2}\ang x\anf K^{-2}.
$$
We will need the following facts on the 3D-calculus $(\Gamma, \dd)$ 
(see \cite{W} and \cite[Subsections 14.1.3 and 14.3.3]{KS}). 
The quantum tangent space of $(\Gamma, \dd)$ is spanned by the functionals
$$
X_0 = q^{-1/2} FK,\quad X_1 = (1 - q^{-2})^{-1} (1 - K^4),\quad 
X_2 = q^{1/2} EK. 
$$
The functionals $f^i_j$ are $f^0_0 = f^2_2 = K^2$, $f^1_1 = K^4$ 
and zero otherwise. The basis elements $\omega_0$, $\omega_1$, 
$\omega_2$ of $\Gamma_{\inv}$ satisfy the relations
$$
\omega^2_0 = \omega^2_1  = \omega^2_2  = 0, \ \ 
\omega_1\omega_0 = -q^4 \omega_0 \omega_1, \ \ 
\omega_2\omega_0 = -q^2 \omega_0\omega_2, \ \ 
\omega_2\omega_1 = -q^4 \omega_1 \omega_2.
$$
The corresponding higher order calculus satisfies the assumptions 
stated at the beginning of this section with $n = 3$ and the volume form 
$\omega := \omega_0 \omega_1 \omega_2$. Further, we have
$$
\omega a = q^{-4} a\omega, \quad
\omega b = q^4 b\omega, \quad
\omega c = q^{-4} \omega, \quad
\omega d = q^4 d \omega
$$
so that $\omega x = \sigma_1 (x) \omega$  with 
$\sigma_1 (x) = K^8 \ang x$. Hence the automorphism 
$\sigma = \sigma_2\circ \sigma_1$ is 
$\sigma (x) = K^6 \ang x \anf K^{-2}$, $x\in \cO(\OSU_q(2))$.
From the commutation relations of the basis elements $\omega_i$, we derive 
\begin{align*}
&\pi (\omega_1\omega_0 \omega_2) = \pi (\omega_0 \omega_2 \omega_1) 
= -q^4, \quad \pi (\omega_2 \omega_1 \omega_0) = -q^{10},\\
&\pi (\omega_1 \omega_2 \omega_0)=\pi (\omega_2 \omega_0 \omega_1)=q^6.
\end{align*}
All other elements $\pi(\omega_i \omega_j \omega_k)$ are zero. 
Inserting these facts in (\ref{xdrel}), we obtain the following explicit 
expression for the $\sigma$-twisted cyclic cocycle $\tau_{\omega, h}$:
\begin{align*}
&\tau_{\omega, h} (x_0, x_1, x_2, x_3) 
= h(\pi (x_0 \dd x_1 \wedge \dd x_2 \wedge \dd x_3))\\
&\qquad = h (x_0 (X_0 \ang x_1) [(K^2 X_1 \ang x_2) (K^6 X_2 \ang x_3) 
- q^4 (K^2 X_2 \ang x_2)(K^4 X_1 \ang x_3)]\\
&\qquad\quad + x_0 (X_1 \ang x_1) [q^6 (K^4 X_2 \ang x_2) (K^6 X_0 \ang x_3) 
- q^4 (K^4 X_0 \ang x_2) (K^6 X_2 \ang x_3)]\\
&\qquad\quad\quad
+ x_0 (X_2 \ang x_1) [q^6 (K^2 X_0 \ang x_2) (K^4 X_1 \ang x_3) 
- q^{10} (K^2 X_1 \ang x_2) (K^6 X_0 \ang x_3)] ).
\end{align*}

Now we develop a commutator representation 
$[\cF, \rho(x)]$ of the 3D-calculus. 
We slightly modify the construction from \cite[Section 3]{S1}.
Let $\Hh$ be the Hilbert space completion of $\cO(\OSU_q(2))$ 
with respect to the inner product 
$(x, y)\hspace{-1pt} :=\hspace{-1pt} h(y^\ast x)$, 
$x,y \in \cO(\OSU_q(2))$. There is a $\ast$-re\-pre\-sen\-tation of the 
left cross product $\ast$-al\-ge\-bra
$\cO(\OSU_q(2)) \rti \cU_q(\su_2)$ on the domain $\cO(\OSU_q(2))$ such that 
$x \in \cO(\OSU_q(2))$ acts by left multiplication and $f \in \cU_q(\su_2)$ 
acts by the left action $\ang$. Let $\cC$ be the closure of the 
Casimir operator $C = FE + (q-q^{-1})^{-2} (q K^2 + q^{-1} K^{-2} - 2)$. 
Let $\zeta (z)$ denote the holomorphic function
$$\zeta (z) = \sum\limits^{\infty}_{n=1} n [n/2]^{-2z}_q [n]_q,\quad 
z \in \dC, \ \,\mathrm{Re}\,  z > 1.
$$
\begin{thl} \label{L3}
If $z \in \dC$, $\mathrm{Re}\,  z > 1$, and $x \in \cO(\OSU_q(2))$, 
then the closure of $\cC^{-z} K^2 x$ and $\cC^{-z} K^{-6} x K^8$ 
are of trace class and 
$$
h(x) = \zeta (z)^{-1} \tr\,  \cC^{-z} K^2 x = \zeta (z)^{-1} 
\tr\,  \cC^{-z} K^{-6} x K^8.
$$
\end{thl}
{\bf Proof.}
The first equality is essentially \cite[Theorem 5.7]{SW}. 
Since we use the left crossed product algebra instead of  the 
right crossed product algebra in \cite{SW}, the operator $K^{-2}$ 
in \cite[Theorem 5.7]{SW} has to be replaced by $K^2$. 
The proof of \cite[Theorem 5.7]{SW} assures that $\cC^{-z} K^2 x$ 
is of trace class for all $x \in \cO(\OSU_q(2))$. 
Given an  $x \in \cO(\OSU_q(2))$,  
there is a $y \in \cO(\OSU_q(2))$ such that $xK^8 = K^8 y$. 
Hence $\cC^{-z} K^{-6} x K^8 = \cC^{-z} K^2 y$ is of trace class. 
The second equality follows by a slight modification of the proof of 
\cite[Theorem 5.7]{SW}. \hfill $\Box$

\sn
Let $\gamma_q$, $\eta_0$, $\eta_1$, $\eta_2$ be complex square matrices 
of the same dimension such that 

\begin{equation}\label{trpro}
\tr\,  \gamma_q \eta_i \eta_j \eta_k = \pi (\omega_i \omega_j \omega_k),\quad 
i, j, k \in \{0, 1, 2\}, 
\end{equation}
and $\eta^\ast_0 = \eta_2$, $\eta^\ast_1 = \eta_1$. 
For instance, take
$$
\gamma_q = \left( \begin{matrix} 1 &0 &0 &0\\
                                 0 &1 &0 &0\\
                                 0 &0 &q^{6} &0\\
                                 0 &0 &0 &q^{6}     \end{matrix} \right),~~
 \eta_0 = \left( \begin{matrix}  0 &0 &0 &1\\
                                 0 &0 &1 &0\\
                                 0 &0 &0 &0\\
                                 0 &0 &0 &0     \end{matrix} \right),~~
\eta_1 = \left( \begin{matrix}   \alpha_1 &0 &0 &0\\
                                 0 &\alpha_2 &0 &0\\
                                 0 &0 &\alpha_3 &0\\
                                 0 &0 &0 &\alpha_4    \end{matrix} \right),~~
$$
and $\eta_2 = \eta^\ast_0$, 
where $\alpha_1$, $\alpha_2$, $\alpha_3$, $\alpha_4$ are real numbers 
such that
$$
\alpha_3 + \alpha_4 = 1,\quad  
\alpha_1 + \alpha_2 = -q^4, \quad 
q^6 (\alpha^3_1 + \alpha^3_2) + \alpha^3_3 + \alpha^3_4 = 0.
$$
Let $\rho$ denote the $\ast$-re\-pre\-sen\-tation of $\cO(\OSU_q(2))$ 
on $\Hh \otimes \dC^4$ given by $\rho (x) := x \otimes I$. 
For simplicity, we write $\cC^{-z}$ for $\cC^{-z} \otimes I$,
$K$ for $K \otimes I$ and $\gamma_q$ for $I\otimes \gamma_q$. 
Since $X^\ast_0 = X_2$, $X^\ast_1 = X_1$ and 
$\eta^\ast_0 = \eta_2$, $\eta^\ast_1 = \eta_1$, the operator 
$\cF := \sum^2_{j=0} X_j \otimes \eta_j$ is self-adjoint on the Hilbert space 
$\Hh \otimes \dC^4$. Set $\Omega_k := \sum^2_{j=0} f^k_j \otimes \eta_j$. 
From (\ref{xfcomma}), we obtain for $x \in \cO(\OSU_q(2))$
\begin{align*}
\Omega_k \rho (x) &= \sum_j f^k_j x \otimes \eta_j 
= \sum_{j,l} (f^k_l \ang x) f^l_j \otimes \eta_j 
= \sum_l (f^k_l \ang x) \Omega_l,\\
[\cF, \rho (x)] &= \sum_k [X_{k}, x] \otimes \eta_k 
= \sum_{k,j} (X_j \ang x) f^j_k \otimes \eta_k = \sum_j (X_j \ang x) \Omega_j.
\end{align*}
Hence, by (\ref{dxfcom}), there is an injective linear mapping 
$\cI : \Gamma \rightarrow \cL(\cO(\OSU_q(2)) \otimes \dC^4)$ 
such that $\cI (x \dd y) = \rho (x) \im [\cF, \rho (x)]$ for all 
$x, y \in \cO(\OSU_q(2))$. This means that the pair $(\cF, \rho)$ 
gives a faithful commutator representation \cite{S1} of the 3D-calculus. 
Further, we compute
\begin{align*}
&\tr_{\Hh \otimes \dC^4}\, \cC^{-z} K^{-6} \gamma_q \rho (x_0) 
[\cF, \rho (x_1)] [\cF, \rho (x_2)] [\cF, \rho (x_3)]\\
&= \tr_{\Hh\otimes \dC^4}\, \cC^{-z} K^{-6} 
\Big( \sum_{i,j,k} x_0 [X_i, x_1] [X_j, x_2] [X_k, x_3] 
\otimes \gamma_q \eta_i \eta_j \eta_k \Big) \\
&= \tr_{\Hh}\, \cC^{-z} K^{-6} 
\Big( \sum_{i,j,k} x_0 [X_i, x_1] [X_j, x_2] [X_{k},x_3] 
\pi (\omega_i \omega_j \omega_k)\Big) \\
&= \tr_{\Hh}\, \cC^{-z} K^{-6} 
\Big(\sum_{i,j,k,l,r,s} x_0 
[X_i, x_1] [X_j, x_2] [X_k,x_3] S(f^{(i,j,k)}_{(l,r,s)}) 
\pi (\omega_l \omega_r \omega_s)  \Big) K^8\\
&= \tr_{\Hh}\, \cC^{-z} K^{-6} 
\pi (x_0 \dd x_1 \wedge \dd x_2 \wedge \dd x_3) K^8\\
&= \zeta (z) h( \pi (x_0 \dd x_1 \wedge \dd x_2 \wedge \dd x_3)) 
= \zeta (z) \tau_{\omega,h} (x_0, x_1, x_2, x_3).
\end{align*}
Here the first equality is the definition of $\cF$ and the second 
equality follows from (\ref{trpro}). For the third equality we 
used the fact that $f^{(i,j,k)}_{(l,r,s)} = K^8$ for all indices 
$i$, $j$, $k$, $l$, $r$, $s$ $\in \{0, 1, 2\}$ for which 
$\pi (\omega_l \omega_r \omega_s) \ne 0$. 
(This is certainly not true for other calculi). 
The fourth equality is the assertion of Lemma \ref{T3/1}, 
while the fifth follows from Lemma \ref{L3}. 
Remind that $\pi (x_0 \dd x_1 \wedge \dd x_2 \wedge \dd x_3)$ 
belongs to $\cO(\OSU_q(2))$. 
The last equality is the definition of the cocycle $\tau_{\omega, h}$. 

Summarizing the preceding, we have proved the following
\begin{tht} \label{T4}
Retaining  the foregoing notation, the $\sigma$-twisted cyclic cocycle 
$\tau_{\omega,h}$ associated with the volume form $\omega$ 
of the 3D-calculus on the Hopf $\ast$-al\-ge\-bra $\cO(\OSU_q(2))$ is given by
\begin{multline*}
\tau_{\omega,h} (x_0, x_1, x_2, x_3) \\
= \zeta (z)^{-1} \tr_{\Hh\otimes \dC^4}\,  \cC^{-z} K^{-6} \gamma_q 
\rho (x_0) [\cF, \rho (x_1)] [\cF, \rho (x_2)] [\cF, \rho (x_3)],
\end{multline*}
where $x_0$, $x_1$, $x_2$, $x_3 \in \cO(\OSU_q(2))$ and $z \in \dC$, 
$\mathrm{Re}\, z>1$.
\end{tht}

\mn

\end{document}